\UseRawInputEncoding
\documentclass[11pt,twoside]{article}
\usepackage{amsmath, amsthm, amsfonts, amssymb,mathtools}
\usepackage{mathrsfs}
\allowdisplaybreaks

\newfam\msbfam
\font\tenmsb=msbm5    \textfont\msbfam=\tenmsb \font\sevenmsb=msbm5
\scriptfont\msbfam=\sevenmsb \font\fivemsb=msbm5
\scriptscriptfont\msbfam=\fivemsb

\newfam\bigfam
\font\tenbig=msbm5 scaled \magstep2   \textfont\bigfam=\tenbig
\font\sevenbig=msbm7 scaled \magstep2 \scriptfont\bigfam=\sevenbig
\font\fivebig=msbm5 scaled \magstep2
\scriptscriptfont\bigfam=\fivebig

\def\dsum{\displaystyle\sum}

\newtheorem{thm}{Theorem}[section]
\newtheorem{lem}[thm]{Lemma}

\newtheorem{rem}[thm]{Remark}
\newtheorem{cor}[thm]{Corollary}
\newtheorem{defn}[thm]{Definition}

\makeatletter

\begin{document}
	\title{\bf 	Lower bound estimates for the rank of universal quadratic forms in some families of real cubic fields with density one}
	
	\author{\bf Liwen Gao and Xuejun Guo}
	
	\date{}
	\maketitle
\renewcommand{\thefootnote}{}
\footnotetext{The authors are supported by National	Nature Science Foundation of China (Nos. 11971226, 12231009).}

	\begin{minipage}{13.5cm}
		
	 {\quad ABSTRACT. In this paper, we establish the explicit lower bound estimates for the rank of universal quadratic forms in some certain families of real cubic fields under the condition of density one. The more general results that represent all multiples of a given rational integer are obtained for totally positive definite quadratic lattices. Our main tools are some properties of indecomposable integers with trace in these fields and short vectors in quadratic lattices. }
	 
\end{minipage}

	{\centering\section{\hspace{-0.3cm}{\bf} Introduction \label{s1}}}
	​A well-known result in number theory is Lagrange's four-square theorem, which states that every natural number can be represented as a sum of four integer squares. To extend the theorem to all positive definite quadratic form, Dickson gave the definition that positive definite quadratic form $f$ over $\mathbb{Z}$ is called universal if $f$ represents all positive integers. A surprisingly simple algorithm was given by Conway and Schneeberger and simplified by Bhargava in \cite{b}, where it was shown that a classical positive quadratic form is universal if and only if it represents all positive integers up to $15$. 
	
	A natural generalization is from $\mathbb{Z}$  to the ring of integers in a general number field.
	Ross proved that there is no ternary positive definite universal form over $\mathbb{Z}$ in \cite{r} and later Earnest and Khosravani gave a more conceptual proof in \cite[Lemma 3]{ek}. However, this result does not hold for any real number fields. In \cite{m}, Maass found that $x^2+y^2+z^2$ is a positive definite universal quadratic form over the ring of integers of $ \mathbb{Q}(\sqrt{5})$.
	Siegel~\cite{s1} in~1945 showed that any diagonsal quadratic form is universal only over the number fields $\mathbb{Q}$ and $\mathbb{Q}(\sqrt 5)$.
	Chan, Kim and Raghavan further identified all possible quadratic real fields admitting a positive definite universal ternary quadratic form and listed all positive definite universal ternary quadratic forms in \cite{ckr}. An analogy of Conway and Schneeberger's 15-theorem was established by Lee in \cite{l} for $ \mathbb{Q}(\sqrt{5})$. Kim and Earnest \cite{ek} showed that there are only finitely many classes of positive definite quadratic forms over a given totally real field with the minimal rank by means of the reduction theory. There are many other articles concerning these subjects such as \cite{hhx,o,xz}.

	Let $K$ be a number field and $O_{K}$ the ring of integers. To study the universal problem in $O_{K}$ conveniently, we need to put it in the more general setting, that is in lattice. A lattice is \textit{universal} if it represents all the totally positive elements of $O_{K}$.
	(see Section 2 for the more precise definition)
	
	There are various interesting problems and conjectures about the universal lattices. For example,  it is well known that  there are only finitely many number fields admitting a ternary universal lattice conjectured by Kitaoka. The conjecture remains unsolved until now, although some studies have been made in \cite{ek,ky2}.
	
	Hsia, Kitaoka, and Kneser in 1978 established a version of Local-global principle for universal forms over number fields in \cite{hkk}. In particular, their Theorem 3 implies that a (totally positive) universal form exists over every totally real number field $K$. The natural problems is to determine the minimal rank for existing positive definite universal quadratic forms, which leads to the following definition of rank of real number field naturally,
	
	\begin{defn}
	 As a universal classical lattice exists over every $H$, we can denote by $R_{class}(H)$ its smallest possible rank over $H$ (and by $R(H)$ the smallest rank without the classical assumption).
\end{defn}
	
	The works of Blomer and Kala \cite{bk1,bk2,k} brought a breakthrough in the area by establishing that the ranks $R(\mathbb{Q}(\sqrt D))$ can be arbitrarily large. Their results were based on the observation that it is hard for a quadratic lattice to represent certain \textit{indecomposable} elements, and so their number can be used to obtain lower bounds for $R(\mathbb{Q}(\sqrt D))$.
	
	To overcome these difficulties, Kala, Yaysyna and $\dot{Z}$mijs gave a lower bound on the rank of the universal quadratic lattice of all real quadratic fields with density one in \cite{kyz}, which they proved that for almost all squarefree $D>0$ and for any $\varepsilon >0$, $R_{class}(\mathbb{Q}(\sqrt D))\geq D^{\frac{1}{12} -\varepsilon}$ and $R(\mathbb{Q}(\sqrt D))\geq D^{\frac{1}{24}-\varepsilon}.$

	As there are a fruit of results for the rank of all real quadratic fields, we can consider whether extend these results to the real cubic field. Through some researches, we find that the proof of results in \cite{kyz} need the behavior of the indecomposables with trace as the key point. Meanwhile, it turned out to appear exactly that the behavior of the indecomposables with trace in some families of real cubic fields in \cite{kt}. We denote 
	all integers that are more than $6$ by $\mathbb{Z}_{\geq 7}$ and the cubic field 
	that is generated by $\rho$ and $1$ by $\mathbb{Q}(\rho)$ where $\rho$ is the root of some cubic polynomial of integer coefficients.
	
	Motivated by the above disscusions,  in this paper we will give lower bound for the rank of the universal quadratic lattices of some families of real cubic fields with density one.
	Now we recall the definitions of two typical families of real cubic fields. 
	In 1974, Shank gave a following family of cubic field in \cite{s2}.
		\begin{defn}
	Shanks' family of the \textit{simplest cubic fields} is defined as $K=\mathbb{Q}(\rho)$, where $\rho$ is the largest root of the polynomial $$x^3-ax^2-(a+3)x-1\text{ for some fixed }a\in \mathbb{Z}_{\geq 7}.$$
		\end{defn}
	
	These fields have a number of advantages: they are Galois, possess units of all signatures, and the ring of integers  equals  for a positive density of $a$, namely, at least for those $a$ for which the square root of the discriminant $\Delta^{1/2}=a^2+3a+9$ is squarefree.
	
	In 2004, Ennola gave a following family of cubic field that is possible to express the fundamental units explicitly and to control the indecomposables of these fields in \cite{e1}.
		\begin{defn}
	Ennola's cubic field is defined as $K=\mathbb{Q}(\rho)$ with the parameter $a\in \mathbb{Z}_{\geq 7}$, where $\rho$ is the smallest root of the polynomial $$x^3+(a-1)x^2-ax-1.$$
\end{defn}

	Our main result of universal form about two families of cubic fields is the following theorems,
	
	\begin{thm} \label{1.4}
		Let $\epsilon>0$. For almost all $a>7$, we have that
		$$R_{class}(\mathbb{Q}(\rho))\geq a^{2-2\epsilon}\quad\text{and}\quad R(\mathbb{Q}(\rho))\geq a^{2-2\epsilon}.$$
	\end{thm}
	By ``almost all'' we mean that such $a$ have (natural) density 1 among the set of all integer $a>7$. 
	
	\begin{thm} \label{1.5}
		Let $\epsilon>0$. For almost all $a>7$ and Ennola's cubic field $\mathbb{Q}(\rho)$, we have that
		$$R_{class}(\mathbb{Q}(\rho))\geq a^{2-2\epsilon}\quad\text{and}\quad R(\mathbb{Q}(\rho))\geq a^{2-2\epsilon}.$$
		``almost all'' is as Theorem 1.4.
	\end{thm}

	This paper will be organized as follows. In Section 2, we will list the necessary definitions and lemmas. In Section 3, we will give the low bound of rank for the Shank's family cubic field with density $1$. In Section 4, the estimates with density $1$ for the rank of Ennola's family will be given. Finally, it will turn out to exist another family that have the similar lower bound for the rank with density $1$ in Section 5.

	{\centering\section{\hspace{-0.3cm}{\bf} Preliminaries \label{s2}}}
	Let $V$ be an $n$-dimensional vector space over $\mathbb{Q}$ endowed with a non-degenerate symmetric bilinear form
	$$\mathcal{B}:V\times V\rightarrow \mathbb{Q},\ \ \text{with the associated quadratic form} \ \  Q(v)=\mathcal{B}(v,v),\quad x\in V .$$
	
	A quadratic $\mathbb{Z}$-lattice $\varLambda\subset V$ is a $\mathbb{Z}$-submodule such that $\mathbb{Q}\varLambda=V$; $\varLambda$ is classical if $\mathcal{B}(v,w)\in \mathbb{Z}$ for all $v,w\in\varLambda$, and positive definite if $Q(v)>0$ for all $0\neq v\in \varLambda$.
	For every quadratic form $Q$ over $\mathbb{Q}$ there is a corresponding quadratic $\mathbb{Z}$-lattice $(\mathbb{Z}^{r},Q)$, where $r$ is the number of variables of $Q$. The associated {\textit Gram matrix} is defined as $G:=(\mathcal{B}(e_{i},e_{j}))$, where 
	$e_{i}$ are elementary vectors in $\mathbb{R}^{r}$. $\Gamma(x)$ denotes the usual gamma function.
	
	Meanwhile, we can extend these definitions about $\mathbb{Z}$-lattice to the cubic field with the ring of integers. Let 
	$H=\mathbb{Q}(\rho)$ and $\mathcal{O}_{H}$ be the integer ring of $H$. Let $V$ be an $n$-dimensional vector space over $H$ equipped with asymmetric bilinear form $$\mathcal{B}:V\times V\rightarrow H,\ \ \text{with the associated quadratic form} \ \  Q(v)=\mathcal{B}(v,v), \quad x\in V .$$
	for $v\in V$. A quadratic $\mathcal{O}_{H}$-lattice $\varLambda\subset V$ is $\mathcal{O}_{H}$-submodule such that $H\varLambda=V$; $\varLambda$ is classical if $\mathcal{B}(v,w)\in \mathcal{O}_{H}$ for all $v,w\in\varLambda$, and positive definite if $Q(v)>0$ for all $0\neq v\in \varLambda$.
	
	We say that $\alpha$ is totally positive if for any real embedding $\sigma:H\rightarrow C$
	satisfying $\sigma(\alpha)>0$ and denote this by $\alpha\succ 0$. All totally positive elements of $\mathcal{O}_{H}$ are denoted by $\mathcal{O}_{H}^{+}$.
	$(\varLambda,Q)$ is totally positive definite if $Q(v)\succ 0$ for every non-zero $v\in \varLambda$. The codifferent of $\mathcal{O}_{H}$ is $\mathcal{O}_{H}^{\vee}=\{\alpha\in H |\mathrm{Tr}_{H/Q}(\alpha \mathcal{O}_{H})\subset\mathbb{Z}\}$, where $\mathrm{Tr}_{H/Q}$: $H\rightarrow Q$ is the trace map. We denote all totally positive elements of $\mathcal{O}_{H}^{\vee}$ by $\mathcal{O}_{H}^{\vee,+}$.
	
	Suppose that we study only totally positive definite quadratic lattices, we still denote them by $\mathcal{O}_{H}$-lattices. A lattice $\varLambda$ represents $\alpha\in\mathcal{O}_{H}^{+}$ if $Q(v)=a$ for some $v\in \varLambda$. If it represents all elements of $\mathcal{O}_{H}^{+}$, then it is universal. Furthermore, a quadratic $\mathcal{O}_{H}$-lattice is $k\mathcal{O}_{H}$-uniserval if it representes all elements of 
	$k\mathcal{O}_{H}^{+}$, for some fixed positive integer $k$.
	
	We say that a totally positive element $\alpha$ of $\mathcal{O}_{H}$ is the indecomposable element if it cannot be written as the sum of $\alpha=\gamma+\delta$ where $\gamma,\delta\in\mathcal{O}_{H}^{+}$. We will use some properties of the indecomposable elements as our main tool in this paper. 

	In \cite{kyz}, the useful estimate of short vectors in $\mathbb{Z}$-lattice is the another  important tool in our paper,

\begin{lem}(see \cite{kyz}). \label{2.1}
	Let $D$ be a classical positive definite $\mathbb{Z}$-lattice of rank $r$ with the Gram matrix $G$ of the basis of $D$ and $N(n)$ be the number of elements $v\in D$ such that $Q(v)=n$. Then $$N(n)\leq C(r,n).$$
	where $C(r,n)$ is defined by
	\begin{equation*}
		C(r,n):=\begin{dcases}
			2r&\text{if } n=1,\\
			\max\{480,2r(r-1)\}& \text{if } n=2,\\
			\frac{\pi^{\frac{r}{2}}}{\Gamma(\frac{r}{2}+1)}\frac{n^{\frac{r}{2}}}{\sqrt{\det(G)}}+ \dsum_{m=0}^{r-1}\binom{r}{m}\frac{\pi^{\frac{m}{2}}}{\Gamma(\frac{m}{2}+1)}n^{\frac{m}{2}}&\text{if } n\geq 3.
		\end{dcases}
	\end{equation*}
	For further use, we denote
	$$B_{1}(r,n):=\frac{1}{2}C(3r,n),$$
	$$B_{2}(r,n):=\frac{1}{2}C(3r,2n).$$
\end{lem}
	
	{\centering\section{\hspace{-0.3cm}{\bf} Shank's cubic fields \label{s3}}}
	By \cite{kt}, we can know the indecomposable elements of the simplest cubic field are almost the form $-v-w\rho+(v+1)\rho^2$ where $0\leq v\leq a$ and $v(a+2)+1\leq w\leq (v+1)(a+1)$.  We need the totally positive elements of these indecompositions to get the main result of the paper and get the following lemma by the proof of \cite[Lemma 5.1]{kt}
	
	\begin{lem}\label{3.1}
		Let $a\in \mathbb{Z}_{\geq 7}$ and $0\leq v\leq a$, $v(a+2)+1\leq w\leq (a+1)(v+1)$.
		Then the elements of the form $-v-w\rho+(v+1)\rho^{2}$ are totally positive.
	\end{lem}
	
	\begin{rem}\label{3.2}
		Under the above assumption of the index scale of $v,w$, we can count the number of these elements
		$$a'=\frac{(a+1)(a+2)}{2}.$$
	\end{rem}
	
	Meanwhile, by the simple calculation, we have the following observation:
	
	\begin{lem}\label{3.3}
		For every $B\geq 1$, we have that if $X>\max\{B,7\}$, then
		$$\#\{\enspace7\leq a\leq X\enspace|\enspace a'\leq B\}<\sqrt{2X}.$$
	\end{lem}
	
	We are now ready to prove the main theorem of this section,
	
	\begin{thm}\label{3.4}
		Let $K$ be the simplest cubic field with parameter $a\in \mathbb{Z}_{\geq 7}$ such that 
		$\mathcal{O}_{K}=\mathbb{Z}(\rho)$ and $a'$ be as in Remark 2.3.
		Assume that $\varLambda$ is an $k\mathcal{O}_{K}$-universal classical quadratic $\mathcal{O}_{K}$-lattice.
		Then $$a'<B_{1}(R,n).$$
	\end{thm}
	
	\noindent{{\it Proof.}} Let $A(v,w)=-v-w\rho+(v+1)\rho^{2}$ where
	$0\leq v\leq \frac{a-1}{3}$, $v(a+3)+1\leq w\leq a(v+1)$. By Lemma 3.1 and \cite[Lemma 5.1]{kt}, each $A(v,w)\succ 0$ and there is $\delta\in \mathcal{O}_{K}^{\vee,+}$ such that $Tr(A(v,w)\delta)=1$.
	
	Fix an integral basis for $\mathcal{O}_{K}$, we can view $\varLambda$ as a classical $\mathbb{Z}$-lattice of rank $3R$ via an isomorphism $\phi:\varLambda\rightarrow \mathbb{Z}^{3R}$. Then we can obtain the classical and positive definite quadratic form $q(v):=Tr(\delta Q(\phi^{-1}(v)))$ over $\mathcal{O}_{K}$.
	
	By our assumption, $Q$ represents all of $k\mathcal{O}_{K}^{+}$. Especially, we can find the vectors $w_{r}$ such that $Q(w_{r})=kA(v,w)$. Set 
	$v_{r}=\phi(w_{r})$  and we have $q(\pm v_{r})=q(v_{r})=Tr(\delta Q(w_{r}))=Tr(k\delta A(v,w))=k$. So each pair of $\pm v_{r}$ are the vectors of norm $k$ of the quadratic lattice $\mathbb{Z}^{3R}$.
	Therefore, the number $N(k)$ of vectors of norm $k$ in our lattice $(\mathbb{Z}^{3R},q)$ satisfies $N(k)\geq 2a'.$ By Lemma 2.2, the desired inequality is obtained. $\hfill\square$

	\begin{thm}\label{3.5}
		Let $R,k$ be positive integers and $X\in \mathbb{Z}_{\geq 7}$ and denote
		$$\mathcal{N}(R,k,X):=\#\{7\leq a\leq X|\exists\enspace k\mathcal{O}_{K} \text{-universal classical lattice of rank $R$ in $K$}\}.$$
		Let $B:=B_{1}(R,n)$ be the quantity defined in Lemma 2.1. Then for all $X>\max\{B,7\}$, we have 
		$$\mathcal{N}(R,k,X)<\sqrt{2X}.$$
	\end{thm}

	\noindent{{\it Proof.}} Theorem 3.4 gives that $a'<B_{1}(R,n)$ , and the result then follows from Lemma 3.3.$\hfill\square$
	
	\begin{cor}\label{3.6}
		For every $\epsilon>0$ and sufficientlt large X, we have 
		$$\#\{\enspace7\leq a\leq X\enspace|\enspace R_{class}(\mathbb{Q}(\rho))\leq a^{2-2\epsilon}\}<\sqrt{2}X^{1-\epsilon},$$
		$$\#\{\enspace7\leq a\leq X\enspace|\enspace R(\mathbb{Q}(\rho))\leq a^{2-2\epsilon}\}<\sqrt{2}X^{1-\epsilon}.$$
	\end{cor}
	
	\medskip
	\noindent{{\it Proof.}}
	The first estimate can be induced by Theorem 3.5.
	
	For the second estimate, if there is a non-classical universal lattice $(\varLambda,Q)$ over $\mathbb{Q}(\rho)$, then $(\varLambda,2Q)$ is $\mathcal{O}_{\mathbb{Q}(\rho)}$-lattice of the same rank and by Theorem 3.5. So we can obtain the desired result.$\hfill\square$
	
	{\centering\section{\hspace{-0.3cm}{\bf} Ennola's cubic fields \label{s4}}}
	Next, we will focus on the Ennola's cubic field. Firstly, we give the totally positive indecomposables of the form $1+w\rho+\rho^{2}$ of Ennola's cubic field in the following lemma by \cite[Proposition 8.1]{kt}.

	\begin{lem}\label{4.1}
		Let $a\in \mathbb{Z}_{\geq 7}$ and $1\leq w\leq a-1$. Then the elements of the form $1+w\rho+\rho^{2}$ are totally positive.
	\end{lem} 
	
	\begin{rem}\label{4.2}
		
		Under the above assumption of the index scale of $w$, we can count the number of these elements
		$$a'=a.$$
	\end{rem} 
	Meanwhile, by the simple calculation, we have the following observation:
	
	\begin{lem}\label{4.3}
		For every $B\geq 1$, we have that if $X>\max\{B^{2},7\}$, then
		$$\#\{\enspace7\leq a\leq X\enspace|\enspace a'\leq B\}<\sqrt{X}.$$
	\end{lem}
	
	We are now ready to prove the main theorem of this section,
	
	\begin{thm}\label{4.4}
		Let $K$ be the Ennola's cubic field with parameter $a\in \mathbb{Z}_{\geq 7}$ such that 
		$\mathcal{O}_{K}=\mathbb{Z}(\rho)$ and $a'$ be as in Remark 4.2.
		Assume that $\varLambda$ is an $m\mathcal{O}_{K}$-universal classical quadratic $\mathcal{O}_{K}$-lattice.
		Then $$a'<B_{2}(R,n).$$
	\end{thm} 
	
	\noindent{{\it Proof.}} Let $C(w)=1+w\rho+\rho^{2}$ where $3\leq w\leq a-1$. By Lemma 4.1 and \cite[Section 8.1]{kt}, each $C(w)\succ 0$ and there is $\delta\in \mathcal{O}_{K}^{\vee,+}$ such that $Tr(C(w)\delta)=2$.
	
	Fix an integral basis for $\mathcal{O}_{K}$, we can view $\varLambda$ as a classical $\mathbb{Z}$-lattice of rank $3R$ via an isomorphism $\phi:\varLambda\rightarrow \mathbb{Z}^{3R}$. Then we can obtain the classical and positive definite quadratic form $q(v):=Tr(\delta Q(\phi^{-1}(v)))$ over $\mathcal{O}_{K}$.
	
	By our assumption, $Q$ represents all of $k\mathcal{O}_{K}^{+}$. Especially, we can find the vectors $w_{r}$ such that $Q(w_{r})=kC(w)$. Set 
	$v_{r}=\phi(w_{r})$  and we have $q(\pm v_{r})=q(v_{r})=Tr(\delta Q(w_{r}))=Tr(k\delta C(w))=2k$. So each pair of $\pm v_{r}$ are the vectors of norm $2k$ of the quadratic lattice $\mathbb{Z}^{3R}$.
	Therefore, the number $N(2k)$ of vectors of norm $2k$ in our lattice $(\mathbb{Z}^{3R},q)$ satisfies $N(2k)\geq 2a'.$ By Lemma 2.2, the desired inequality is obtained. $\hfill\square$
	
	\begin{thm}\label{4.5}
		Let $R,k$ be positive integers and $X\in \mathbb{Z}_{\geq 7}$ and
		denote
		$$\mathcal{N}(R,k,X):=\#\{\enspace7\leq a\leq X|\enspace\exists\enspace k\mathcal{O}_{K} \text{-universal classical lattice of rank $R$ in $K$ } \}.$$
		
		Let $B:=B_{2}(R,n)$ be the quantity defined in Lemma 2.2. Then for all $X>\max\{B,7\}$, we have 
		$$\mathcal{N}(R,k,X)<\sqrt{X}.$$
	\end{thm}
	
	\noindent{{\it Proof.}} Theorem 4.4 gives that $a'<B_{2}(R,n)$ , and the result then follow from Lemma 4.3.$\hfill\square$
	
	\begin{cor}\label{4.6}
		For every $\epsilon>0$, sufficientlt large X and Ennola's cubic field $\mathbb{Q}(\rho)$, we have 
		$$\#\{\enspace7\leq a\leq X\enspace|\enspace R_{class}(\mathbb{Q}(\rho))\leq a^{2-2\epsilon}\}<X^{1-\epsilon},$$
		$$\#\{\enspace7\leq a\leq X\enspace|\enspace R(\mathbb{Q}(\rho))\leq a^{2-2\epsilon}\}<X^{1-\epsilon}.$$
	\end{cor}  
	
	\noindent{{\it Proof.}}
	The first estimate can be induced by Theorem 4.5.
	
	For the second estimate, if there is a non-classical universal lattice $(\varLambda,Q)$ over $\mathbb{Q}(\rho)$, then $(\varLambda,2Q)$ is $\mathcal{O}_{\mathbb{Q}(\rho)}$-lattice of the same rank and by Theorem 4.5. So we can get the desired result.$\hfill\square$.
	
	{\centering\section{\hspace{-0.3cm}{\bf} Another family of cubic field \label{s5}}}
	We can get another family of cubic field by letting $r=a$ and $s=a+2$ in the latter family \cite{t}. The definition is as follow,
	
	\begin{defn}\label{5.1}
		The field is defined as $K=\mathbb{Q}_(\rho)$ with the parameter $a\in \mathbb{Z}_{\geq 7}$, where $\rho$ is the smallest root of the polynomial $$x^3-(2a+2)x^2+a(a+2)x-1.$$
	\end{defn}
	
	Next, we give the totally positive indecomposables of the form $-1+((a+2)w+1)\rho-w\rho^{2}$ of the above field by \cite[Proposition 8.3]{kt}.
	
	\begin{lem}\label{5.2}
		Let $a\in \mathbb{Z}_{\geq 7}$ and $a\leq w\leq 2a-1$. Then the elements of the form $-1+((a+2)w+1)\rho-w\rho^{2}$ are totally positive.
	\end{lem}
	
	\begin{rem}\label{5.3}
		Under the above assumption of the index scale of $w$, we can count the number of these elements
		$$a'=a.$$
	\end{rem}  
	Meanwhile, by the simple calculation, we have the following observation:
	
	\begin{lem}\label{5.4}
		For every $B\geq 1$, we have that if $X>\max\{B^{2}+1,7\}$, then
		$$\#\{\enspace7\leq a\leq X\enspace|\enspace a'\leq B\}<\sqrt{X}.$$
	\end{lem}   
	We are now ready to prove the main theorem of this section,
	
	\begin{thm}\label{5.5}
		Let $K$ be the cubic field given by Definition 5.1 with parameter $a\in \mathbb{Z}_{\geq 7}$ such that 
		$\mathcal{O}_{K}=\mathbb{Z}(\rho)$ and $a'$ be as in Remark 5.3.
		Assume that $\varLambda$ is an $k\mathcal{O}_{K}$-universal classical quadratic $\mathcal{O}_{K}$-lattice.
		Then $$a'<B_{2}(R,n).$$
	\end{thm}  
	
	\noindent{{\it Proof.}} 
	Let $D(w)=-1+((a+2)w+1)\rho-w\rho^{2}$ where $a\leq w\leq 2a-1$. By lemma 5.2 and \cite[Section 8.2]{kt}, each $D(w)\succ 0$ and there is $\delta\in \mathcal{O}_{K}^{\vee,+}$ such that $Tr(D(w)\delta)=2$.
	
	Fix an integral basis for $\mathcal{O}_{K}$, we can view $\varLambda$ as a classical $\mathbb{Z}$-lattice of rank $3R$ via an isomorphism $\phi:\varLambda\rightarrow \mathbb{Z}^{3R}$. Then we can obtain the classical and positive definite quadratic form $q(v):=Tr(\delta Q(\phi^{-1}(v)))$ over $\mathcal{O}_{K}$.
	
	By our assumption, $Q$ represents all of $k\mathcal{O}_{K}^{+}$. Especially, we can find the vectors $w_{r}$ such that $Q(w_{r})=kD(w)$. Set 
	$v_{r}=\phi(w_{r})$  and we have $q(\pm v_{r})=q(v_{r})=Tr(\delta Q(w_{r}))=Tr(k\delta D(w))=2k$. So each pair of $\pm v_{r}$ are the vectors of norm $2k$ of the quadratic lattice $\mathbb{Z}^{3R}$.
	Therefore, the number $N(2k)$ of vectors of norm $2k$ in our lattice $(\mathbb{Z}^{3R},q)$ satisfies $N(2k)\geq 2a'.$ By Lemma 2.2, the desired inequality is obtained. $\hfill\square$

	\begin{thm}\label{5.6}
		Let $R,k$ be positive integers and $X\in \mathbb{Z}_{\geq 7}$ and denote
		$$\mathcal{N}(R,k,X):=\#\{\enspace7\leq a\leq X|\enspace\exists\enspace k\mathcal{O}_{K} \text{-universal classical lattice of rank $R$ in $K$ } \}.$$
		Let $B:=B_{2}(R,n)$ be the quantity defined in Lemma 2.2. Then for all $X>\max\{B^{2}+1,7\}$, we have 
		$$\mathcal{N}(R,k,X)<\sqrt{X}.$$
	\end{thm}
	
	\noindent{{\it Proof.}} Theorem 5.5 gives that $a'<B_{2}(R,n)$ , and the result then follow from Lemma 5.4.$\hfill\square$
	
	\begin{cor}\label{5.7}
		For every $\epsilon>0$, sufficientlt large X and Ennola's cubic field $Q(\rho)$, we have 
		$$\#\{\enspace7\leq a\leq X\enspace|\enspace R_{class}(\mathbb{Q}(\rho))\leq a^{2-2\epsilon}\}<X^{1-\epsilon},$$
		$$\#\{\enspace7\leq a\leq X\enspace|\enspace R(\mathbb{Q}(\rho))\leq a^{2-2\epsilon}\}<X^{1-\epsilon}.$$
	\end{cor}
	
	\noindent{{\it Proof.}}
	The first estimate can be induced by Theorem 5.6.
	
	For the second estimate, if there is a non-classical universal lattice $(\varLambda,Q)$ over $\mathbb{Q}(\rho)$, then $(\varLambda,2Q)$ is $\mathcal{O}_{\mathbb{Q}(\rho)}$-lattice of the same rank and by Theorem 5.6. So we can get the desired result.$\hfill\square$
	
	\begin{cor}\label{5.8}
		Let $\epsilon>0$. For almost all $a>7$ and the cubic field $\mathbb{Q}(\rho)$ given by definition 5.1 we have that
		$$R_{class}(\mathbb{Q}(\rho))\geq a^{2-2\epsilon}\quad\text{and}\quad R(\mathbb{Q}(\rho))\geq a^{2-2\epsilon}.$$
		``almost all'' is the same as Theorem 1.4.
	\end{cor}
	
{99}

	\medskip
	
	Liwen Gao\par
	Department of Mathematics\par
	Nanjing University\par
	Nanjing 210093, CHINA \par
	Email address: {\bf gaoliwen1206@smail.nju.edu.cn}\par
	
	\centerline{}
	Xuejun Guo\par
	Department of Mathematics\par
	Nanjing University\par
	Nanjing 210093, CHINA \par
	Email address: {\bf guoxj@nju.edu.cn}\par
	
\end{document}